\documentclass[thmsa,a4,12pt]{article}
\usepackage{amscd}
\usepackage{amsmath}
\usepackage{graphicx}
\usepackage{amsfonts}
\usepackage{amssymb}
\newtheorem{theorem}{Theorem}

\newtheorem{corollary}[theorem]{Corollary}

\newtheorem{definition}[theorem]{Definition}

\newtheorem{lemma}[theorem]{Lemma}

\newtheorem{proposition}[theorem]{Proposition}

\begin{document}

\title{Approximation of recurrence in negatively curved metric spaces }
\author{Charalambos Charitos and Georgios Tsapogas}
\date{December 22, 1998}
\maketitle
\begin{abstract}
For metric spaces with curvature less than or equal to $\chi,$ $\chi<0$, it is
shown that a recurrent geodesic is approximated by closed geodesics. A counter
example is provided for the converse.
\end{abstract}

\section{Introduction and preliminaries}

In hyperbolic geometry it has been shown \footnotetext{1991
\textit{Mathematics Subject Classification.} Primary 53C22; Secondary
53C23.}lately that many geometric properties are determined by the distance
function on the space itself rather than the differential structure. It is
shown in this work that, partially, this is the case with the notion of
recurrence. For complete hyperbolic manifolds, a recent result of Aebisher,
Hong and McCullough (see \cite{A-H-M}) states that a geodesic is recurrent if
and only if it is approximated by closed geodesics. We show that, in metric
spaces with curvature less than or equal to $\chi,$ $\chi<0$, recurrent
geodesics are approximated by closed geodesics (see theorem \ref{main} below).
The proof of the converse statement crucially depends on the manifold
structure, in particular on the fact that two geodesics coincide if they do so
on an open interval. Hence, the converse statement fails in our context due to
the bifurcation property of geodesics. A counter example exhibiting this
failure is provided in section \ref{ce} below. A geodesic $\gamma$ is called
\textit{recurrent} if there exists a sequence $\left\{  t_{n}\right\}
\subset\mathbb{R},t_{n}\rightarrow\infty$ such that $t_{n}\gamma
\rightarrow\gamma$ as $t_{n}\rightarrow\infty.$ Convergence in this definition
is meant to be uniform convergence on compact sets which, in fact, induces the
topology on the space $GX$ consisting of all (local) isometries $\mathbb{R}%
\rightarrow X$ when $X$ is (not) simply connected. $\mathbb{R}$ acts on $GX$
by right translations, namely, $\left(  t,g\right)  \rightarrow tg,$ where
$tg:\mathbb{R}\rightarrow X$ is the geodesic defined by $tg\left(  s\right)
=g\left(  s+t\right)  ,$ $s\in\mathbb{R}$. This action is simply the geodesic
flow. The notion of convergence in the above definition is analogous to the
tangential condition which defines recurrence in the manifold case. We use the
notion of approximation given in definition \ref{recdef} below which was
introduced in \cite{A-H-M} in order to characterize recurrent geodesics in
hyperbolic manifolds.

$X$ will always denote a locally compact, complete, geodesic metric space with
curvature less than or equal to $\chi,$ $\chi<0.$ Recall that a geodesic
metric space is said to have curvature less than or equal to $\chi$ if each
$x\in X$ has a neighborhood $V_{x}$ such that every geodesic triangle of
perimeter strictly less than $\frac{2\pi}{\sqrt{\chi}}$ (=+$\infty$ when
$\chi\leq0$) contained in $V_{x}$ satisfies $CAT-\left(  \chi\right)  $
inequality (see \cite{Pau} for definitions and basic properties). We will
denote the metric by $d\left(  \cdot,\cdot\right)  $ and will use the same
letter to denote distance when the metric space to which we refer is
understood. All curves are assumed to be parametrized by arclength. A
\textit{geodesic segment} in $X$ is an isometry $c:I\rightarrow X$, where $I$
is a closed interval in $\mathbb{R}$. A \textit{geodesic} in $X $ is a map
$c:\mathbb{R}\rightarrow X$ such that for each closed interval $I\subset
\mathbb{R}$, the map $c\,|_{I}:I\rightarrow X$ is a geodesic segment. A
\textit{local geodesic segment }(usually called \textit{geodesic arc}) in X is
a map $c:I\rightarrow X$ such that for each $t\in I$ there is an
$\varepsilon>0$ such that $c\,|_{\left[  t-\varepsilon,t+\varepsilon\right]
\cap I}:\left[  t-\varepsilon,t+\varepsilon\right]  \cap I\rightarrow X$ is a
geodesic segment. Similarly, a \textit{local geodesic }$\mathbb{R}\rightarrow
X$ is defined. A \textit{closed geodesic} in X is a local geodesic
$c:\mathbb{R}\rightarrow X$ which is a periodic map.

\begin{definition}
\label{appdef}An oriented geodesic $g$ in $X$ is said to be approximated by
closed geodesics if, for every $\varepsilon>0$ and every $x\in$%
\textrm{Im\thinspace}$g$, there exists a closed oriented geodesic $c$ such
that for some point $y\in$\textrm{Im\thinspace}$c$,
\[
d\left(  c\left(  t+t_{y}\right)  ,g\left(  t+t_{x}\right)  \right)
<\varepsilon
\]
for all $t\in\left[  0,period\left(  c\right)  \right]  ,$ where $t_{x}%
,t_{y}\in\mathbb{R}$ with $x=g\left(  t_{x}\right)  $ and $y=c\left(
t_{y}\right)  .$
\end{definition}

The following theorem is the main result of this paper.

\begin{theorem}
\label{main}Let $X$ be a locally compact, complete, geodesic metric space
which has curvature less than or equal to $\chi,$ $\chi<0.$ If a geodesic or
geodesic ray in $X$ is recurrent then it is approximated by closed geodesics.
\end{theorem}

The proof of theorem \ref{main} uses the notion of quasi-geodesic and its
stability properties. We will closely follow notation and terminology
appearing in \cite[Ch.3]{C-D-P} where we refer the reader for first
definitions and basic properties of quasi-geodesics. Here we only recall the
following definition.

\begin{definition}
\label{quasig}Let $f:\left[  a,b\right]  \rightarrow X$ be a continuous map
with $-\infty\leq a\leq b\leq+\infty$ and $\lambda,\kappa,L$ real numbers with
$\lambda\geq1,\kappa\geq0,L>0$. \newline $f$ is a $\left(  \lambda
,\kappa,L\right)  -$quasi-geodesic if for every subinterval $\left[
a^{\prime},b^{\prime}\right]  $ of $\left[  a,b\right]  $ satisfying
\[
lengthf\left(  \left[  a^{\prime},b^{\prime}\right]  \right)  \leq L,
\]
the following inequality holds
\[
length\,\,f\left(  \left[  a^{\prime},b^{\prime}\right]  \right)  \leq\lambda
d\left(  f\left(  a^{\prime}\right)  ,f\left(  b^{\prime}\right)  \right)
+\kappa
\]
\end{definition}

The next proposition is a well know fact for $CAT-\left(  \chi\right)  $
spaces. We include a short proof of it, since it is difficult to find exact
reference (when $X$ is a geometric polyhedron this result follows from
\cite[p.403]{Bri})

\begin{proposition}
\label{localg}Let $M$ be a complete geodesic space satisfying $CAT-\left(
\chi\right)  $ inequality with $\chi<0.$ Every local geodesic segment in $M$
is a geodesic segment.
\end{proposition}

\noindent\textbf{Proof.} Let $\delta:\left[  0,L\right]  \rightarrow M$ be a
local geodesic segment in $M$, $L>0$. Set
\[
\mathit{l}=\sup\left\{  t\in\left[  0,L\right]  \bigm\vert\delta\,|_{\left[
0,t\right]  }\,\,is\,\,a\,\,geodesic\,\,segment\right\}
\]
Apparently, $\mathit{l}>0$ and by completeness of $M,$ $\delta\,|_{\left[
0,\mathit{l}\right]  }$ is a geodesic segment joining $\delta\left(  0\right)
$ with $\delta\left(  \mathit{l}\right)  $. Assuming the conclusion is not
true, i.e. $\mathit{l}<L$, let $\varepsilon$ be a positive number such that
$\delta\,|_{\left[  \mathit{l}-\varepsilon,\mathit{l}+\varepsilon\right]  }$
is a geodesic segment. Denote by $\left[  \delta\left(  0\right)
,\delta\left(  \mathit{l}+\varepsilon\right)  \right]  $ the geodesic segment
in $M$ joining $\delta\left(  0\right)  $ with $\delta\left(  \mathit{l}%
+\varepsilon\right)  $. Since $\delta\,|_{\left[  0,\mathit{l}+\varepsilon
\right]  }$ is not the geodesic segment joining $\delta\left(  0\right)  $
with $\delta\left(  \mathit{l}+\varepsilon\right)  ,$
\begin{equation}
d\left(  \delta\left(  0\right)  ,\delta\left(  \mathit{l}+\varepsilon\right)
\right)  <d\left(  \delta\left(  0\right)  ,\delta\left(  \mathit{l}\right)
\right)  +d\left(  \delta\left(  \mathit{l}\right)  ,\delta\left(
\mathit{l}+\varepsilon\right)  \right)  \label{pai}%
\end{equation}
The points $\delta\left(  0\right)  ,\delta\left(  \mathit{l}\right)  $ and
$\delta\left(  \mathit{l}+\varepsilon\right)  $ define a geodesic triangle in
$M $. Denote by $\Delta=\left(  \overline{\delta\left(  0\right)  }%
,\overline{\delta\left(  \mathit{l}\right)  },\overline{\delta\left(
\mathit{l}+\varepsilon\right)  }\right)  $ the corresponding comparison
triangle which is non-degenerate by inequality (\ref{pai}). Choose points $B$
on $\delta\,|_{\left[  0,\mathit{l}\right]  }$ and $B^{\prime}$ on
$\delta\,|_{\left[  \mathit{l},\mathit{l}+\varepsilon\right]  }$ such that
$d\left(  B,\delta\left(  \mathit{l}\right)  \right)  =d\left(  B^{\prime
},\delta\left(  \mathit{l}\right)  \right)  =\varepsilon^{\prime}<\varepsilon$
and denote by $\overline{B}$ and $\overline{B^{\prime}}$ the corresponding
points on the comparison triangle. Then by (\ref{pai}) the angle of $\Delta$
at $\overline{\delta\left(  \mathit{l}\right)  }$ is smaller that $\pi$ and
therefore
\[
d\left(  \overline{B},\overline{B^{\prime}}\right)  <d\left(  \overline
{B},\overline{\delta\left(  \mathit{l}\right)  }\right)  +d\left(
\overline{\delta\left(  \mathit{l}\right)  },\overline{B^{\prime}}\right)
=2\varepsilon^{\prime}%
\]
By comparison, $d\left(  B,B^{\prime}\right)  \leq d\left(  \overline
{B},\overline{B^{\prime}}\right)  $ so we obtain
\[
d\left(  B,B^{\prime}\right)  <d\left(  B,\delta\left(  \mathit{l}\right)
\right)  +d\left(  \delta\left(  \mathit{l}\right)  ,B^{\prime}\right)
\]
This contradicts the fact that $\delta\,|_{\left[  \mathit{l}-\varepsilon
^{\prime},\mathit{l}+\varepsilon^{\prime}\right]  }$ is a geodesic segment.
\hfill\rule{1mm}{3mm}

Let $\widetilde{X}$ be the universal cover of $X$ and $p:\widetilde
{X}\rightarrow X$ the projection map. $\widetilde{X}$ becomes a metric space
as follows : given $\widetilde{x},\widetilde{y}\in\widetilde{X}$ choose any
curve $\widetilde{c}:\left[  a,b\right]  \rightarrow\widetilde{X}$ with
$\widetilde{c}\left(  a\right)  =\widetilde{x}$ and $\widetilde{c}\left(
b\right)  =\widetilde{y}$ and define the distance from $\widetilde{x}$ to
$\widetilde{y}$ to be the length of the unique length minimizing curve in the
homotopy class of $p\widetilde{c}$ with endpoints fixed. For the existence of
the length minimizing curve see \cite{Gro2}. This distance function is a
metric on $\widetilde{X}$ which inherits the properties of $X,$ namely,
$\widetilde{X}$ becomes a complete geodesic locally compact (hence, proper)
metric space. $\pi_{1}\left(  X\right)  $ acts on $\widetilde{X}$ and the
action commutes with $p.$ As the projection $p$ is a local isometry, it
follows that $\pi_{1}\left(  X\right)  $ acts on $\widetilde{X}$ by local
isometries. Using the fact that $\widetilde{X}$ is geodesic and proposition
\ref{localg}, it is routine to show that $\pi_{1}\left(  X\right)  $ acts on
$\widetilde{X}$ by isometries. In addition, $\widetilde{X}$ has curvature less
than or equal to $\chi,$ $\chi<0$ and, by a theorem of Gromov (see for example
\cite[p.325]{Pau}), $\widetilde{X}$ satisfies $CAT-\left(  \chi\right)  $ inequality.

$GX$ is by definition the space of all local geodesics $\mathbb{R}\rightarrow
X$ and, by proposition \ref{localg} above, $G\widetilde{X}$ is the space
consisting of all global geodesics $\mathbb{R}\rightarrow\widetilde{X}.$ The
topology on these spaces is uniform convergence on compact sets. The boundary
$\partial\widetilde{X}$ can be defined using either equivalence classes of
sequences or, equivalence classes of geodesic rays. The local compactness
assumption on $X$ implies that $\widetilde{X}$ is proper and hence the two
definitions coincide (see \cite[Ch.2]{C-D-P}). We will be using them
interchangeably. For any two distinct points $\xi,\eta$ in $\partial
\widetilde{X}$ there exists a unique, up to parametrization, (oriented)
geodesic $g$ with $g\left(  -\infty\right)  =\xi$ and $g\left(  \infty\right)
=\eta$ (see for example \cite[prop. 2]{Cha}). We need the following lemma
which asserts that the projection of a point onto a geodesic always exists.

\begin{lemma}
\label{pro}Let $g$ be a geodesic in $G\widetilde{X}$ (or a geodesic segment)
and $x_{0}$ a point in $\widetilde{X}.$ There exists a unique real number $s$
such that $g\left(  s\right)  $ realizes the distance of $x_{0}$ from $\left.
\mathrm{Im\,}g,\right.  $ i.e. $dist\left(  x_{0},\mathrm{Im\,}g\right)
=d\left(  x_{0},g\left(  s\right)  \right)  .$
\end{lemma}

\noindent\textbf{Proof }We may assume that $x_{0}\notin$\textrm{Im\thinspace
}$g.$ Existence is apparent. Assume that $s\neq s^{\prime}$ are two such
numbers. The points $g\left(  s\right)  ,g\left(  s^{\prime}\right)  $ and
$x_{0} $ define a non-degenerate geodesic triangle in $\widetilde{X}$ and
denote by $\Delta=\left(  \overline{g\left(  s\right)  },\overline{g\left(
s^{\prime}\right)  },\overline{x_{0}}\right)  $ the corresponding comparison
triangle. $\Delta$ is an equilateral triangle in the unique complete simply
connected Riemannian 2-manifold of constant sectional curvature $\chi.$ Hence,
the angles of $\Delta$ at $\overline{g\left(  s\right)  }$ and $\overline
{g\left(  s^{\prime}\right)  }$ are each less than $\pi/2.$ Therefore, there
exists a point $\overline{g\left(  t\right)  }$ on the side of $\Delta$
opposite to $x_{0}$ such that $d\left(  x_{0},\overline{g\left(  t\right)
}\right)  <d\left(  x_{0},\overline{g\left(  s\right)  }\right)  =d\left(
x_{0},\overline{g\left(  s^{\prime}\right)  }\right)  .$ By $CAT-\left(
\chi\right)  $ inequality, $d\left(  x_{0},g\left(  t\right)  \right)  \leq
d\left(  x_{0},g\left(  s\right)  \right)  ,$ a contradiction.\hfill\rule{1mm}{3mm}

\noindent\textbf{Remark 1 }If $c\in GX$ is a closed geodesic and $x_{0}\in X,$
the same argument applied to a lifting $\widetilde{c}$ of $c$ shows that there
exists a unique point $B\in$\textrm{Im}$\,c$ such that $d\left(
x_{0},B\right)  =dist\left(  x_{0},\mathrm{Im}\,c\right)  .$

\noindent\textbf{Remark 2 }Set $\partial^{2}\widetilde{X}=\left\{  \left(
\xi,\eta\right)  \in\partial\widetilde{X}\times\partial\widetilde{X}:\xi
\neq\eta\right\}  $ and let $\rho:G\widetilde{X}\rightarrow\partial
^{2}\widetilde{X}$ be the fiber bundle given by $\rho\left(  g\right)
=\left(  g\left(  -\infty\right)  ,g\left(  +\infty\right)  \right)  $. Since
for any two distinct points $\xi,\eta$ in $\partial\widetilde{X}$ there exists
a unique (oriented) geodesic $g$ with $g\left(  -\infty\right)  =\xi$ and
$g\left(  \infty\right)  =\eta$ (see for example \cite[prop. 2]{Cha}), the
fiber of $\rho$ is $\mathbb{R}.$ Moreover, this bundle is trivial (see for
example \cite[Th. 4.8]{Cham}). To define a trivialization, let $x_{0}$ be
a base point and let
\begin{equation}
H:G\widetilde{X}\overset{\approx}{\longrightarrow}\partial^{2}\widetilde
{X}\times\mathbb{R}\label{lization}%
\end{equation}
be the trivialization of $\rho$ with respect to $x_{0}$ defined by
\[
H\left(  g\right)  =\left(  g\left(  -\infty\right)  ,g\left(  +\infty\right)
,s\right)
\]
where $-s$ is the real number provided by lemma \ref{pro}. Note that the
composite of the geodesic flow $\mathbb{R}\times G\widetilde{X}\rightarrow
G\widetilde{X}$ with $H$ is given by the formula
\[
\left(  \xi_{1},\xi_{2},s\right)  \longrightarrow\left(  \xi_{1},\xi
_{2},s+t\right)
\]
$\,$for all $\left(  \xi_{1},\xi_{2}\right)  \in\partial^{2}\widetilde{X}$ and
$s\in\mathbb{R}.$

\section{Recurrent geodesics}

\begin{definition}
\label{recdef}A geodesic $\gamma$ in $X$ is called recurrent if there exists a
sequence $\left\{  t_{n}\right\}  \subset\mathbb{R},t_{n}\rightarrow\infty$
such that $t_{n}\gamma\rightarrow\gamma$ as $t_{n}\rightarrow\infty.$
\end{definition}

For a recurrent geodesic $\gamma$ in $X$ there exists \textit{a sequence of
closed} (in fact, piece-wise geodesic) \textit{curves} $\left\{  \gamma
_{n}\right\}  _{n\in\mathbb{N}}$, \textit{associated to} $\gamma$ as follows :
fix a convex neighborhood $U$ of $\gamma\left(  0\right)  ,$ i.e. a
neighborhood which satisfies the following property : for all $x,y\in U$ there
exists unique geodesic segment with endpoints $x$ and $y$ lying entirely in
$U.$ Such a neighborhood exists (see for example \cite{Bal}). If $\left\{
t_{n}\right\}  $ is the sequence given by definition \ref{recdef} above and
$\varepsilon_{n}=d\left(  \gamma\left(  0\right)  ,\gamma\left(  t_{n}\right)
\right)  $, let $K\in\mathbb{N}$ such that $\gamma\left(  t_{n}\right)  \in U$
for all $n\geq K.$ Define $\gamma_{n},n\geq K$ to be the curve
\begin{equation}
\gamma_{n}:\left[  0,t_{n}+\varepsilon_{n}\right]  \rightarrow X \label{assoc}%
\end{equation}
with $\gamma_{n}\left(  t\right)  =\gamma\left(  t\right)  \,\,\forall
\,t\in\left[  0,t_{n}\right]  $ and $\gamma_{n}|_{\left[  t_{n},t_{n}%
+\varepsilon_{n}\right]  }$ the unique geodesic segment in $U$ joining
$\gamma\left(  t_{n}\right)  $ with $\gamma\left(  0\right)  .$ Note that
$t_{n}+\varepsilon_{n}$ is the period of the closed curve $\gamma_{n}$. In the
sequel, we will refer to these closed curves by writing $\gamma_{n}%
,n\in\mathbb{N}$ but it will always be implicit that $n$ is large enough so
that $\gamma_{n}$ are defined.

Using the following lemma, we may assume that given a recurrent geodesic
$\gamma$, the associated closed curves $\left\{  \gamma_{n}\right\}
_{n\in\mathbb{N}}$ are not homotopic to a point.

\begin{lemma}
\label{loxodro}Given a recurrent geodesic $\gamma$ there exists $M\in
\mathbb{N}$ such that each closed curve $\gamma_{n}$, $n\in\mathbb{N}$
associated to $\gamma$ is not homotopic to a point, provided $n\geq M.$
\end{lemma}

\noindent\textbf{Proof.} Let $\widetilde{\gamma}$ be a lift of $\gamma$ to the
universal cover $\widetilde{X}$ of $X$ parametrized so that $\widetilde
{\gamma}\left(  0\right)  $ projects to $\gamma\left(  0\right)  =\gamma
_{n}\left(  0\right)  .$ The curve $\gamma_{n}|_{\left[  0,t_{n}\right]  }$ is
a local geodesic segment and, by proposition \ref{localg}, its lift
$\widetilde{\gamma_{n}}|_{\left[  0,t_{n}\right]  }$ to $\widetilde{X}$
starting at $\widetilde{\gamma}\left(  0\right)  $ is a geodesic segment.
Moreover, $\gamma_{n}|_{\left[  t_{n},t_{n}+\varepsilon_{n}\right]  }$ and its
lift $\widetilde{\gamma_{n}}|_{\left[  t_{n},t_{n}+\varepsilon_{n}\right]  }$
to $\widetilde{X}$ starting at $\widetilde{\gamma}\left(  0\right)  $ are both
geodesic segments. We have
\[%
\begin{array}
[c]{lll}%
d\left(  \widetilde{\gamma}\left(  t_{n}+\varepsilon_{n}\right)
,\widetilde{\gamma}\left(  0\right)  \right)   & \geq &  d\left(
\widetilde{\gamma}\left(  t_{n}\right)  ,\widetilde{\gamma}\left(  0\right)
\right)  -d\left(  \widetilde{\gamma}\left(  t_{n}+\varepsilon_{n}\right)
,\widetilde{\gamma}\left(  t_{n}\right)  \right)  \\
& = & t_{n}-\varepsilon_{n}\,
\end{array}
\]
Since $\varepsilon_{n}\rightarrow0$ and $t_{n}\rightarrow\infty$ as
$n\rightarrow\infty$ we may choose $M\in\mathbb{N}$ such that $\widetilde
{\gamma}\left(  t_{n}+\varepsilon_{n}\right)  ,\widetilde{\gamma}\left(
0\right)  $ are distinct for all $n\geq M.$ Therefore, $\widetilde{\gamma_{n}%
}|_{\left[  0,t_{n}+\varepsilon_{n}\right]  },$ which is the lift of the
closed curve $\gamma_{n}$ starting at $\widetilde{\gamma}\left(  0\right)
=\widetilde{\gamma_{n}}\left(  0\right)  ,$ has distinct endpoints and,
therefore, $\gamma_{n},n\geq M$ is not homotopic to a point. \hfill\rule{1mm}{3mm}

The following proposition shows that the lifts (to the universal cover
$\widetilde{X}$) of the closed curves $\gamma_{n}$ associated to a recurrent
geodesic $\gamma$ are, for $n$ large enough, quasi-geodesics with arbitrarily
large $L$. Recall that a\ $CAT-\left(  \chi\right)  $ space is a $\delta
-$hyperbolic space in the sense of Gromov (see for example \cite[Sec.2]{Pau}).
This applies to the universal covering $\widetilde{X}$, since it satisfies
$CAT-\left(  \chi\right)  $ inequality globally. Let $\delta$ denote the
hyperbolicity constant of the space $\widetilde{X}.$

\begin{proposition}
\label{bound}Let $\gamma$ be a recurrent geodesic in $X$ and $\left\{
\gamma_{n}\right\}  _{n\in\mathbb{N}}$ the associated closed curves. For every
$L>0,$ there exists $N\in\mathbb{N}$ such that all lifts $\widetilde
{\gamma_{n}}:\mathbb{R}\rightarrow\widetilde{X}$ of $\gamma_{n}$ with $n\geq
N$ are $\left(  \lambda,\kappa,L\right)  -$quasi-geodesics provided
$\kappa>16\delta$ and $\lambda=1,$ where $\delta$ is the hyperbolicity
constant of $\widetilde{X}.$
\end{proposition}

\noindent\textbf{Proof.} Let $\gamma$ be a recurrent geodesic and $L>0$ be
given. The sequence $\left\{  t_{n}\right\}  $ given by definition
\ref{recdef} converges to infinity. Moreover, $\varepsilon_{n}=d\left(
\gamma\left(  0\right)  ,\gamma\left(  t_{n}\right)  \right)  \rightarrow0$
and $t_{n}+\varepsilon_{n}=period\left(  \gamma_{n}\right)  $ also converges
to infinity as $n\rightarrow\infty$. Hence, we may choose $N$ such that
\begin{equation}
t_{n}+\varepsilon_{n}>L\,\,and\,\,\varepsilon_{n}<\frac{1}{2}\left(
\kappa-16\delta\right)  \,\,for\,\,all\,\,n\geq N. \label{log}%
\end{equation}
Let now $\left[  a,b\right]  $ be any interval with $b-a<L$ (cf. definition
\ref{quasig}). For each $n\geq N$ there exists an integer $k_{n}$ such that
\begin{equation}
\widetilde{\gamma_{n}}\left(  \left[  a,b\right]  \right)  \subset
\widetilde{\gamma_{n}}\left(  \left.  \left[  \left(  k_{n}-1\right)  \left(
t_{n}+\varepsilon_{n}\right)  ,k_{n}\left(  t_{n}+\varepsilon_{n}\right)
+t_{n}\right]  \right.  \right)  \label{three}%
\end{equation}
Denote by $\left[  \widetilde{\gamma_{n}}\left(  a\right)  ,\widetilde
{\gamma_{n}}\left(  b\right)  \right]  $ the unique geodesic segment in
$\widetilde{X}$ joining $\widetilde{\gamma_{n}}\left(  a\right)  $ with
$\widetilde{\gamma_{n}}\left(  b\right)  $ and set
\[%
\begin{array}
[c]{l}%
A_{n}:=\widetilde{\gamma_{n}}\left(  k_{n}\left(  t_{n}+\varepsilon
_{n}\right)  -\varepsilon_{n}\right) \\
y_{k_{n}}:=\widetilde{\gamma_{n}}\left(  k_{n}\left(  t_{n}+\varepsilon
_{n}\right)  \right)
\end{array}
\]
The distance of any point on $\left[  \widetilde{\gamma_{n}}\left(  a\right)
,\widetilde{\gamma_{n}}\left(  b\right)  \right]  $ form $\widetilde
{\gamma_{n}}\left(  \left[  a,b\right]  \right)  $ is bounded by a number
which depends on the hyperbolicity constant $\delta$ of the space
$\widetilde{X}$and on the number of geodesic segments which constitute
$\widetilde{\gamma_{n}}\left(  \left[  a,b\right]  \right)  ,$ see \cite[Lemma
1.5 p.25]{C-D-P}. In our case here, $\widetilde{\gamma_{n}}\left(  \left[
a,b\right]  \right)  $ consists of at most three geodesic segments (since the
right hand side of inclusion (\ref{three}) above consists of 3 geodesic
segments) and the bound is $8\delta.$ Hence we have
\begin{equation}
d\left(  y_{k_{n}},\left[  \widetilde{\gamma_{n}}\left(  a\right)
,\widetilde{\gamma_{n}}\left(  b\right)  \right]  \right)  \leq8\delta
\label{lemmap25}%
\end{equation}
By lemma \ref{pro}, let $B_{n}$ be the point on $\left[  \widetilde{\gamma
_{n}}\left(  a\right)  ,\widetilde{\gamma_{n}}\left(  b\right)  \right]  $
which realizes the distance in the left hand side of inequality \ref{lemmap25}%
. Assume that neither $\widetilde{\gamma_{n}}\left(  a\right)  $ nor
$\widetilde{\gamma_{n}}\left(  b\right)  $ lies on the geodesic segment
$\left[  A_{n},y_{k_{n}}\right]  .$ Then we have the following triangle
inequalities
\[%
\begin{array}
[c]{l}%
d\left(  \widetilde{\gamma_{n}}\left(  a\right)  ,A_{n}\right)  \leq d\left(
\widetilde{\gamma_{n}}\left(  a\right)  ,B_{n}\right)  +d\left(
B_{n},y_{k_{n}}\right)  +d\left(  y_{k_{n}},A_{n}\right) \\
d\left(  y_{k_{n}},\widetilde{\gamma_{n}}\left(  b\right)  \right)  \leq
d\left(  y_{k_{n}},B_{n}\right)  +d\left(  B_{n},\widetilde{\gamma_{n}}\left(
b\right)  \right)
\end{array}
\]
which, after employing the fact that $d\left(  A_{n},y_{k_{n}}\right)
=\varepsilon_{n},$ become
\[%
\begin{array}
[c]{ll}%
\quad\quad\quad length\,\,\widetilde{\gamma_{n}}\left(  \left[  a,b\right]
\right)  & =d\left(  \widetilde{\gamma_{n}}\left(  a\right)  ,A_{n}\right)
+d\left(  A_{n},y_{k_{n}}\right)  +d\left(  y_{k_{n}},\widetilde{\gamma_{n}%
}\left(  b\right)  \right) \\
& \leq2\varepsilon_{n}+2d\left(  y_{k_{n}},B_{n}\right)  +d\left(
\widetilde{\gamma_{n}}\left(  a\right)  ,\widetilde{\gamma_{n}}\left(
b\right)  \right) \\
_{by\,\,inequality\,\,(\ref{lemmap25})} & \leq2\varepsilon_{n}+2\cdot
8\delta+d\left(  \widetilde{\gamma_{n}}\left(  a\right)  ,\widetilde
{\gamma_{n}}\left(  b\right)  \right) \\
_{by\,\,inequality\,\,(\ref{log})} & \leq\kappa+d\left(  \widetilde{\gamma
_{n}}\left(  a\right)  ,\widetilde{\gamma_{n}}\left(  b\right)  \right)
\end{array}
\]
The case where $\widetilde{\gamma_{n}}\left(  a\right)  $ and/or
$\widetilde{\gamma_{n}}\left(  b\right)  $ lies on $\left[  A_{n},y_{k_{n}%
}\right]  $ is treated similarly. \hfill\rule{1mm}{3mm}

\begin{corollary}
\label{loxo}For $n\in\mathbb{N}$ sufficiently large, the isometry of
$\widetilde{X}$ in $\pi_{1}\left(  X\right)  $ which corresponds to the
homotopy class of the closed curve $\gamma_{n}$ is hyperbolic.
\end{corollary}

\noindent\textbf{Proof.} It suffices to show that each $\widetilde{\gamma_{n}%
}:\mathbb{R}\rightarrow\widetilde{X}$ determines exactly two boundary points
$\widetilde{\gamma_{n}}\left(  -\infty\right)  $, $\widetilde{\gamma_{n}%
}\left(  +\infty\right)  .$ By lemma \ref{bound} there exists an
$M\in\mathbb{N} $ such that $\widetilde{\gamma_{n}}$ is a quasi-isometry for
all $n\geq M.$ Each such $\widetilde{\gamma_{n}}$ induces a map $\partial
\mathbb{R}\rightarrow\partial\widetilde{X}$ which is a homeomorphism onto its
image, see \cite[Th.2.2 p.35]{C-D-P}. As $\partial\mathbb{R}$ consists of two
distinct points, $\widetilde{\gamma_{n}}\left(  -\infty\right)  $,
$\widetilde{\gamma_{n}}\left(  +\infty\right)  \in\partial\widetilde{X}$ are
also distinct for all $n\geq M.$\hfill\rule{1mm}{3mm}

It now follows that a recurrent geodesic $\gamma$ in $X$ as well as each of
the (oriented) closed curves $\gamma_{n},n\geq M$ (cf. lemma \ref{loxodro} and
corollary \ref{loxo} above) determine exactly two boundary points in
$\partial\widetilde{X}$ denoted by $\widetilde{\gamma}\left(  -\infty\right)
$, $\widetilde{\gamma}\left(  +\infty\right)  $ and $\widetilde{\gamma_{n}%
}\left(  -\infty\right)  $, $\widetilde{\gamma_{n}}\left(  +\infty\right)  $
respectively. We need the following lemma concerning these boundary points.
Recall that $\widetilde{X}\cup\partial\widetilde{X}$ is a compact space which
is metrizable (see \cite[p.134]{C-D-P}) and we will denote such metric by
$d_{\widetilde{X}\cup\partial\widetilde{X}}$.

\begin{lemma}
\label{limit}$\widetilde{\gamma_{n}}\left(  -\infty\right)  \rightarrow
\widetilde{\gamma}\left(  -\infty\right)  $ and $\widetilde{\gamma_{n}}\left(
+\infty\right)  \rightarrow\widetilde{\gamma}\left(  +\infty\right)  $ as
$n\rightarrow\infty.$
\end{lemma}

\noindent\textbf{Proof.} As above, let $\varepsilon_{n}=length($%
\textrm{Im\thinspace}$\gamma_{n})-t_{n}$ so that $t_{n}+\varepsilon_{n}$ is
the period of $\gamma_{n}$. We first show that $\widetilde{\gamma_{n}}\left(
+\infty\right)  \rightarrow\widetilde{\gamma}\left(  +\infty\right)  .$
Consider the sequence $\widetilde{\gamma_{n}}\left(  k\left(  t_{n}%
+\varepsilon_{n}\right)  \right)  ,k\in\mathbb{N}$ which converges to
$\widetilde{\gamma_{n}}\left(  +\infty\right)  $ as $k\rightarrow\infty.$
Thus, there exists $k_{n}\in\mathbb{N}$ such that
\begin{equation}
d_{\widetilde{X}\cup\partial\widetilde{X}}\left(  \widetilde{\gamma_{n}%
}\left(  k_{n}\left(  t_{n}+\varepsilon_{n}\right)  \right)  ,\widetilde
{\gamma_{n}}\left(  +\infty\right)  \right)  <1/n. \label{ine}%
\end{equation}
Now consider the sequences $y_{n}:=\widetilde{\gamma_{n}}\left(  k_{n}\left(
t_{n}+\varepsilon_{n}\right)  \right)  $ and $x_{n}:=\widetilde{\gamma}\left(
t_{n}\right)  $, $n\in\mathbb{N}.$ Since $x_{n}\rightarrow\widetilde{\gamma
}\left(  +\infty\right)  $, by inequality (\ref{ine}) above it is enough to
show that the sequences $\left\{  x_{n}\right\}  $ and $\left\{
y_{n}\right\}  $ represent the same element in $\partial\widetilde{X}$ or, in
other words, that the hyperbolic product $\left(  x_{n},y_{n}\right)  _{x_{0}%
}$ with respect to the base point $x_{0}:=\widetilde{\gamma}\left(  0\right)
$ converges to $+\infty$ as $n\rightarrow+\infty.$ For the notion of
hyperbolic product of sequences and their equivalence, see \cite{C-D-P}.

The stability property of quasi-geodesics states (see corollary 1.10 of
\cite[p.31]{C-D-P}) that given any two numbers $\kappa\geq0$ and $\lambda
\geq1,$ there exists a constant $C$ depending on $\lambda,\kappa$ and on the
hyperbolicity constant $\delta$ of the space such that if $L$ is bigger than
$2C$ then every $\left(  \lambda,\kappa,L\right)  -$quasi-geodesic $f:\left[
a,b\right]  \rightarrow\widetilde{X}$ lies within a $C-$neighborhood of the
geodesic segment $\left[  f\left(  a\right)  ,f\left(  b\right)  \right]  .$
By choosing $\lambda=1,$ $\kappa>16\delta$ where $\delta$ is the hyperbolicity
constant of the space $\widetilde{X}$ and $L>2C$ we obtain, by proposition
\ref{bound} above, a natural number $N$ such that all $\widetilde{\gamma_{n}%
}:\mathbb{R}\rightarrow\widetilde{X}$ with $n\geq N$ are $\left(
\lambda,\kappa,L\right)  -$quasi-geodesics. In particular, $\widetilde
{\gamma_{n}}:\left[  0,y_{n}\right]  \rightarrow\widetilde{X}$ are $\left(
\lambda,\kappa,L\right)  -$quasi-geodesics for all $n\geq N.$ Therefore, by
corollary 1.10 of \cite[p.31]{C-D-P} as explained above,
\[
d\left(  x_{n},x_{n}^{\prime}\right)  <C\,\,\,\forall\,\,n\geq N
\]
where $x_{n}^{\prime}$ denotes the projection of $x_{n}$ on the geodesic
segment $\left[  \widetilde{\gamma}\left(  0\right)  ,y_{n}\right]  $ (cf.
lemma \ref{pro}). Hence,
\[%
\begin{array}
[c]{ll}%
\left(  x_{n},y_{n}\right)  _{x_{0}} & =\frac{1}{2}\left(  d\left(
x_{n},x_{0}\right)  +d\left(  y_{n},x_{0}\right)  -d\left(  x_{n}%
,y_{n}\right)  \right) \\
& \geq\frac{1}{2}\left(  d\left(  x_{n}^{\prime},x_{0}\right)  -C+d\left(
y_{n},x_{0}\right)  -d\left(  x_{n}^{\prime},y_{n}\right)  -C\right) \\
& =\left(  x_{n}^{\prime},y_{n}\right)  _{x_{0}}-C\\
& =d\left(  x_{0},x_{n}^{\prime}\right)  -C
\end{array}
\]
Apparently, $d\left(  x_{0},x_{n}^{\prime}\right)  \rightarrow\infty$ as
$n\rightarrow+\infty$ and, hence, $\left(  x_{n},y_{n}\right)  _{x_{0}%
}\rightarrow\infty$ as required.

In order to show that $\widetilde{\gamma_{n}}\left(  -\infty\right)
\rightarrow\widetilde{\gamma}\left(  -\infty\right)  $ we work in a similar
manner: the sequence $\widetilde{\gamma_{n}}\left(  -k\left(  t_{n}%
+\varepsilon_{n}\right)  \right)  ,k\in\mathbb{N}$ converges to $\widetilde
{\gamma_{n}}\left(  -\infty\right)  $ as $k\rightarrow\infty.$ Hence, there
exists $k_{n}\in\mathbb{N}$ such that $d_{\widetilde{X}\cup\partial
\widetilde{X}}\left(  \widetilde{\gamma_{n}}\left(  -k_{n}\left(
t_{n}+\varepsilon_{n}\right)  \right)  ,\widetilde{\gamma_{n}}\left(
-\infty\right)  \right)  <1/n.$ As before, sequences $\left\{  y_{n}\right\}
$ and $\left\{  x_{n}\right\}  $ are defined by $y_{n}:=\widetilde{\gamma_{n}%
}\left(  -k_{n}\left(  t_{n}+\varepsilon_{n}\right)  \right)  $ and
$x_{n}:=\widetilde{\gamma}\left(  -t_{n}\right)  $, $n\in\mathbb{N}.$ Then we
use the same arguments to show that the hyperbolic product $\left(
x_{n},y_{n}\right)  _{x_{0}}$ with respect to the base point $x_{0}%
:=\widetilde{\gamma}\left(  0\right)  $ converges to $+\infty$ as
$n\rightarrow+\infty.$\hfill\rule{1mm}{3mm}

\section{Proof of main theorem}

Let $\gamma$ be a recurrent geodesic, $\varepsilon>0$ and $x\in$%
\textrm{Im\thinspace}$\gamma$ be given. We may assume that $x=\gamma\left(
0\right)  .$ Let $\left\{  t_{n}\right\}  $ be the sequence given by
definition \ref{recdef} and $\left\{  \gamma_{n}\right\}  $ the sequence of
the associated closed curves given by formula (\ref{assoc}) above. For each
$n\in\mathbb{N},$ there exists a unique closed geodesic $c_{n}$ in the free
homotopy class of $\gamma_{n}.$ The number $t_{n}+\varepsilon_{n}$ is the
period of $\gamma_{n}$ and let $s_{n}$ denote the period of $c_{n}$
(apparently, $s_{n}<t_{n}+\varepsilon_{n}$). Let $B_{n}$ be the projection of
$\gamma\left(  0\right)  $ onto \textrm{Im\thinspace}$c_{n},$ i.e. $d\left(
\gamma\left(  0\right)  ,B_{n}\right)  =d\left(  \gamma\left(  0\right)
,\mathrm{Im\,}c_{n}\right)  .$ Such a point exists and is unique by remark 1
following lemma \ref{pro}. Lift $\gamma$ to an isometry $\widetilde{\gamma
}:\mathbb{R}\rightarrow\widetilde{X}$ with a base point $\widetilde{\gamma
}\left(  0\right)  $ satisfying $p\left(  \widetilde{\gamma}\left(  0\right)
\right)  =\gamma\left(  0\right)  ,$ where $p:\widetilde{X}\rightarrow X$ is
the universal covering map. Lift each $c_{n}$ to an isometry $\widetilde
{c_{n}}:\mathbb{R}\rightarrow\widetilde{X}$ and parametrize it so that
$\widetilde{c_{n}}\left(  0\right)  $ is a point $\widetilde{B_{n}}$
satisfying
\[
d\left(  \widetilde{B_{n}},\widetilde{\gamma}\left(  0\right)  \right)
=d\left(  B_{n},\gamma\left(  0\right)  \right)  \,\,and\,\,p\left(
\widetilde{B_{n}}\right)  =B_{n}%
\]
For the reader's convenience, we have gathered all the above notation in
figure 1.

\begin{center}
\begin{figure}[ptb]
\begin{picture}(460,140)
\thinlines
\put(15,20){\line(1,0){385}}
\put(15,100){\line(1,0){410}}
\put(40,100){\line(4,-1){320}}
\put(40,100){\line(6,-1){340}}
\put(342.3,100){\line(2,-3){37.7}}
\put(40,20){\circle*{5}}
\put(240,20){\circle*{5}}
\put(357,20){\circle*{5}}
\put(220,20){\circle*{5}}
\put(43,100){\circle*{5}}
\put(240,100){\circle*{5}}
\put(342.3,100){\circle*{5}}
\put(409,100){\circle*{5}}
\put(230,52.5){\circle*{5}}
\put(380,43.4){\circle*{5}}
\put(235,67.5){\circle*{5}}
\put(0,17){$\widetilde{c_n}$}
\put(3,98){$\widetilde{\gamma}$}
\put(30,110){$\widetilde{\gamma}(0)$}
\put(25,5){$\widetilde{c_n}(0)=\widetilde{B_n}$}
\put(230,110){$\widetilde{\gamma}(s)$}
\put(328.3,110){$\widetilde{\gamma}(t_n)$}
\put(378,110){$\widetilde{\gamma}(t_n+{\epsilon}_n)$}
\put(375,53.4){$\widetilde{{\gamma}_n}(t_n+{\epsilon}_n)$}
\put(235,73){$D_n$}
\put(220,40){$F_n$}
\put(206,5){$F_n^{\prime}$}
\put(235,5){$\widetilde{c_n}(s)$}
\put(339,2){$\widetilde{c_n}(s_n)=\phi_n\left( \widetilde{B_n}\right) $ }
\end{picture}
\caption{{}}%
\label{fig1}%
\end{figure}
\end{center}

\noindent Since
\[
p\left(  \widetilde{\gamma_{n}}\left(  t_{n}+\varepsilon_{n}\right)  \right)
=p\left(  \widetilde{\gamma}\left(  0\right)  \right)  =\gamma\left(
0\right)
\]
and $\gamma_{n},c_{n}$ are homotopic, the isometry $\phi_{n}$ of
$\widetilde{X}$ which translates $\widetilde{c_{n}}$ (in the positive
direction) satisfies
\[
\phi_{n}\left(  \widetilde{\gamma}\left(  0\right)  \right)  =\widetilde
{\gamma_{n}}\left(  t_{n}+\varepsilon_{n}\right)
\]
Moreover,
\begin{equation}%
\begin{array}
[c]{lll}%
d\left(  \widetilde{\gamma_{n}}\left(  t_{n}+\varepsilon_{n}\right)  ,\phi
_{n}\left(  \widetilde{B_{n}}\right)  \right)   & = & d\left(  \phi_{n}\left(
\widetilde{\gamma_{n}}\left(  0\right)  \right)  ,\phi_{n}\left(
\widetilde{B_{n}}\right)  \right)  \\
& = & d\left(  \widetilde{\gamma}\left(  0\right)  ,\widetilde{B_{n}}\right)
\end{array}
\label{last}%
\end{equation}
We now proceed to show that given $\varepsilon>0$, there exists $N\in
\mathbb{N}$ such that for all $n\geq N$%
\begin{equation}
d\left(  \widetilde{\gamma}\left(  s\right)  ,\widetilde{c_{n}}\left(
s\right)  \right)  <\varepsilon\,\,\,\forall\,s\in\left[  0,s_{n}\right]
\label{final}%
\end{equation}
Recall that $s_{n}$ is the period of $c_{n}$ and $s_{n}<t_{n}+\varepsilon
_{n}=period\left(  \gamma_{n}\right)  .$ Using lemma \ref{limit} and the fact
that $\gamma_{n},c_{n}$ are homotopic for all $n$ large enough, we have that
$\widetilde{c_{n}}\left(  +\infty\right)  \rightarrow\widetilde{\gamma}\left(
+\infty\right)  $ and $\widetilde{c_{n}}\left(  -\infty\right)  \rightarrow
\widetilde{\gamma}\left(  -\infty\right).$ Let $H:G\widetilde{X}%
\overset{\approx}{\longrightarrow}\partial^{2}\widetilde{X}\times\mathbb{R}$
be the trivilization of the fiber bundle $G\widetilde{X}\longrightarrow
\partial^{2}\widetilde{X}$ with respect to the base point $x_{0}%
=\widetilde{\gamma}\left(  0\right)  .$ This homeomorphism was described in
remark 2 following lemma \ref{pro}. By the choice of parametrization for each
$\widetilde{c_{n}}$ made above $\bigl($ie $\widetilde{c_{n}}\left(  0\right)
=\widetilde{B_{n}}\bigr)$, we have that $H^{-1}\left(  \widetilde{c_{n}%
}\left(  -\infty\right)  ,\widetilde{c_{n}}\left(  +\infty\right)  ,0\right)
=\widetilde{c_{n}}.$ Moreover, $H^{-1}\left(  \widetilde{\gamma}\left(
-\infty\right)  ,\widetilde{\gamma}\left(  +\infty\right)  ,0\right)
=\widetilde{\gamma}$ and, thus, $\widetilde{c_{n}}\rightarrow\widetilde
{\gamma}$ uniformly on compact sets. Observe that such convergence is weaker
than property (\ref{final}). However, it implies, in particular, that
$dist\left(  \widetilde{\gamma}\left(  0\right)  ,\mathrm{Im\,}\widetilde
{c_{n}}\right)  \rightarrow0$ as $n\rightarrow\infty.$ Hence, we may choose
$N\in\mathbb{N}$ such that
\begin{equation}
d\left(  \widetilde{\gamma}\left(  0\right)  ,\widetilde{B_{n}}\right)
<\varepsilon/5,\,\,for\,\,\,all\,\,\,n\geq N\label{uni}%
\end{equation}
Moreover, we may choose $N$ such that, in addition, the following inequality
is satisfied
\begin{equation}
\varepsilon_{n}=d\left(  \widetilde{\gamma}\left(  t_{n}\right)
,\widetilde{\gamma}\left(  t_{n}+\varepsilon_{n}\right)  \right)
<\varepsilon/5,\,\,for\,\,\,all\,\,\,n\geq N\label{eps}%
\end{equation}
To show inequality (\ref{final}), let $s\in\left[  0,s_{n}\right]  $ be
arbitrary and let $D_{n}$ $\left(  resp.\,\,F_{n}\right)  $ be the point on
the geodesic segment $\left[  \widetilde{\gamma}\left(  0\right)
,\widetilde{\gamma_{n}}\left(  t_{n}+\varepsilon_{n}\right)  \right]  \left(
resp.\,\,\,\left[  \widetilde{\gamma}\left(  0\right)  ,\phi_{n}\left(
\widetilde{B_{n}}\right)  \right]  \right)  $ whose distance from
$\widetilde{\gamma}\left(  0\right)  $ is equal to $s.$ Then,
\[%
\begin{array}
[c]{lll}%
d\left(  \widetilde{\gamma}\left(  s\right)  ,\widetilde{c_{n}}\left(
s\right)  \right)   & \leq &  d\left(  \widetilde{\gamma}\left(  s\right)
,D_{n}\right)  +d\left(  D_{n},F_{n}\right)  +d\left(  F_{n},\widetilde{c_{n}%
}\left(  s\right)  \right)  \\
& \leq &  d\left(  \widetilde{\gamma}\left(  s\right)  ,D_{n}\right)
+d\left(  D_{n},F_{n}\right)  +d\left(  F_{n},F_{n}^{\prime}\right)  +d\left(
F_{n}^{\prime},\widetilde{c_{n}}\left(  s\right)  \right)
\end{array}
\]
where $F_{n}^{\prime}$ is the point on $\left[  \widetilde{B_{n}},\phi
_{n}\left(  \widetilde{B_{n}}\right)  \right]  $ satisfying
\[
d\left(  F_{n},\phi_{n}\left(  \widetilde{B_{n}}\right)  \right)  =d\left(
F_{n}^{\prime},\phi_{n}\left(  \widetilde{B_{n}}\right)  \right)  .
\]
By comparison (see for example \cite[prop.29]{Tr}) we have
\[%
\begin{array}
[c]{ll}%
\begin{array}
[c]{r}%
d\left(  \widetilde{\gamma}\left(  s\right)  ,D_{n}\right)  \\
d\left(  D_{n},F_{n}\right)  \\
d\left(  F_{n},F_{n}^{\prime}\right)  \\
d\left(  F_{n}^{\prime},\widetilde{c_{n}}\left(  s\right)  \right)
\end{array}
&
\begin{array}
[c]{l}%
\leq\,\,d\left(  \widetilde{\gamma}\left(  t_{n}+\varepsilon_{n}\right)
,\widetilde{\gamma_{n}}\left(  t_{n}+\varepsilon_{n}\right)  \right)
\leq\,\,2\varepsilon_{n}\\
\leq\,\,d\left(  \widetilde{\gamma}\left(  t_{n}+\varepsilon_{n}\right)
,\phi_{n}\left(  \widetilde{B_{n}}\right)  \right)  \\
\leq\,\,d\left(  \widetilde{\gamma}\left(  0\right)  ,\widetilde{B_{n}%
}\right)  \\
\leq\,\,\left|  d\left(  \widetilde{\gamma}\left(  0\right)  ,\phi_{n}\left(
\widetilde{B_{n}}\right)  \right)  -d\left(  \widetilde{B_{n}},\phi_{n}\left(
\widetilde{B_{n}}\right)  \right)  \right|  <d\left(  \widetilde{\gamma
}\left(  0\right)  ,\widetilde{B_{n}}\right)
\end{array}
\end{array}
\]
Combining the above inequalities with inequalities (\ref{last}), (\ref{uni})
and (\ref{eps}), we obtain property (\ref{final}) which completes the proof of
the existence of a sequence of closed geodesics approximating a given
recurrent geodesic.\hfill\rule{1mm}{3mm}

\noindent\textbf{Remark }Let $\Gamma$ be a discrete group of isometries of a
locally compact, simply connected, complete geodesic metric space $Y$
satisfying $CAT-\left(  \chi\right)  $ inequality, $\chi<0.$ The notion of
controlled concentration points in the limit set of $\Gamma$ can be defined as
follows. $\xi\in\partial\Gamma$ is a \textit{controlled concentration point}
if it admits a neighborhood $U$ containing $\xi$ with the following property :
for every neighborhood $V$ of $\xi$ there exists an element $\gamma\in\Gamma$
such that $\gamma\left(  U\right)  \subset V$ and $\xi\in\gamma\left(
U\right)  .$ Following \cite{A-H-M}, one can show that $\xi$ is a controlled
concentration point if and only if there exists a sequence of $\left\{
\phi_{n}\right\}  $ of distinct elements of $\Gamma$ such that $\phi
_{n}\left(  \xi\right)  \rightarrow\xi$ and $\phi_{n}\left(  0\right)
\rightarrow\eta$ with $\eta\neq\xi.$ The proof in this more general setting is
identical with the one provided in \cite{A-H-M} except that the convergence
property used there, namely, $\phi_{n}\left(  x\right)  \rightarrow\eta$ for
all $x\in Y\cup\partial Y,$ is provided in our case by proposition 7.2 in
\cite[Ch.1]{Coo1}. The latter property for $\xi$ is equivalent to the
existence of a recurrent geodesic $\gamma$ with $\gamma\left(  +\infty\right)
=\xi$ and $\gamma\left(  -\infty\right)  =\eta.$ Hence we obtain the following
connection between recurrent geodesics and controlled concentration points
which also holds for manifolds (see \cite{A-H-M}).

\begin{theorem}
Let $Y$ be a simply connected, locally compact, complete geodesic metric space
$Y$ satisfying $CAT-\left(  \chi\right)  $ inequality, $\chi<0$ and $\Gamma$ a
discrete group of isometries of $Y.$ A limit point $\xi\in\partial Y$ is a
controlled concentration point if and only if $\gamma\left(  +\infty\right)
=\xi$ for some recurrent geodesic $\gamma$ in $Y.$
\end{theorem}

\section{Construction of a counter example\label{ce}}

\begin{figure}
[ptb]
\begin{center}
\includegraphics[
height=2.8686in,
width=1.8827in
]
{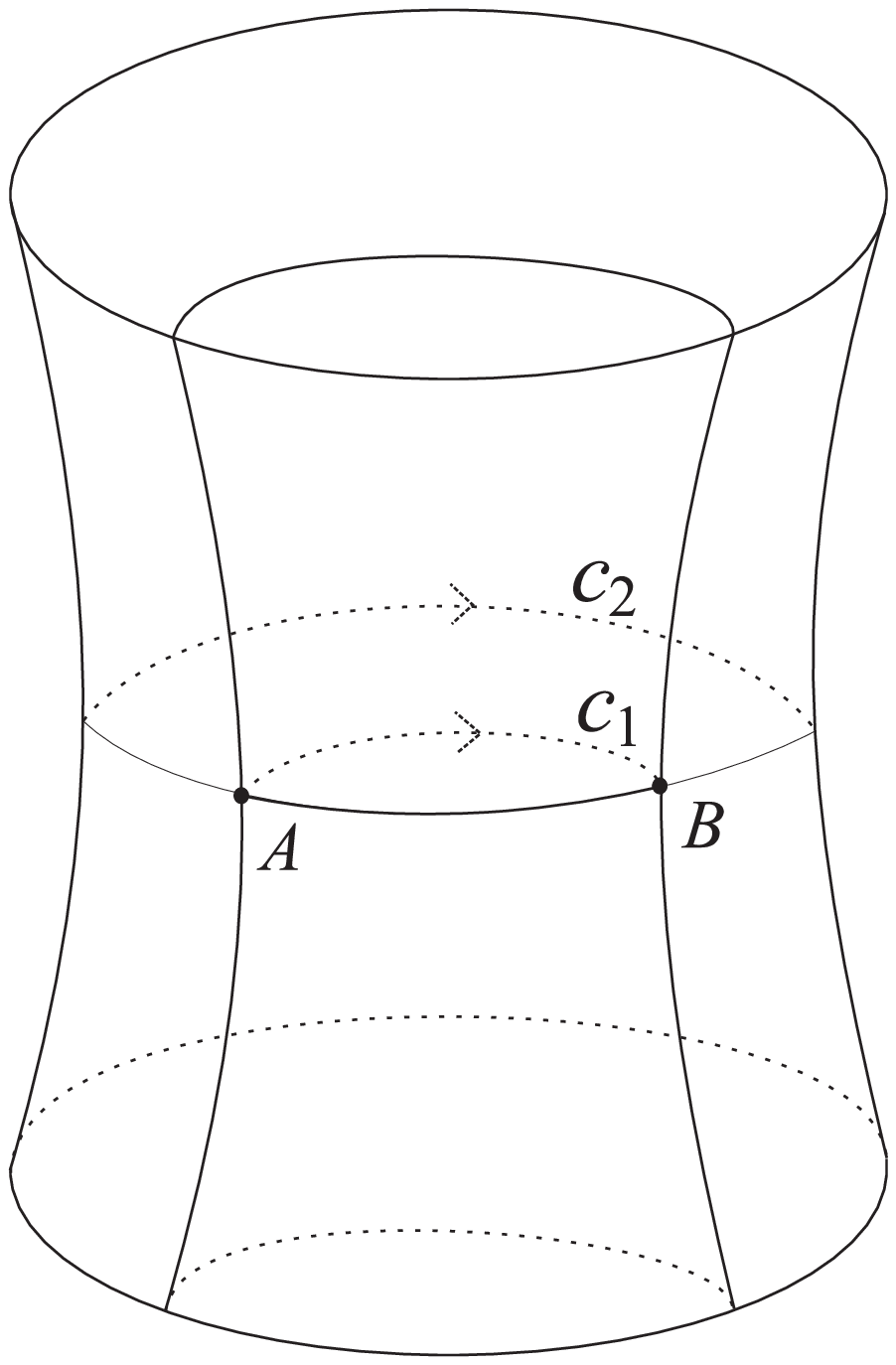}
\label{cylinder}
\end{center}
\end{figure}

As it was mentioned in the introduction approximation by closed geodesics does
not imply recurrence. The following example demonstrates the existence of a
geodesic in a $CAT-\left(  \chi\right)  ,\chi<0$ space which is not recurrent
but it can be approximated by closed geodesics in the sense of definition
\ref{appdef}. Let $X$ be the union of two hyperbolic cylinders identified
along a (convex) geodesic strip bounded by two geodesic segments (see figure
\ref{cylinder}). We may adjust the geometry of $X$ so that the unique simple
closed geodesic in each cylinder, denoted by $c_{1}$ and $c_{2},$ have a
common image in the geodesic strip, namely, the geodesic segment indicated by
letters $A$ and $B$ in figure \ref{cylinder}. Using Cor. 5 of \cite{Bal} and
the fact that the geodesic strip is a convex closed subset it follows that $X$
is a $CAT-\left(  \chi\right)  $ space with $\chi<0. $%

Let $\omega_{1}$ and $\omega_{2}$ be the periods of $c_{1}$ and $c_{2}$
respectively and assume that $c_{1}$ and $c_{2}$ are parametrized so that
$c_{1}\left(  0\right)  =c_{2}\left(  0\right)  =B$ and clockwise i.e.
$c_{1}\left(  s\right)  =c_{2}\left(  s\right)  $ for all $s\in\left[
0,d\left(  A,B\right)  \right]  .$ Define $\gamma:\mathbb{R}\rightarrow X$ as
follows :
\[%
\begin{array}
[c]{l}%
\gamma\left(  t\right)  =c_{1}\left(  t\right)  ,\,\,\,for\,\,\,t\in\left[
0,\omega_{1}\right] \\
\gamma\left(  t\right)  =c_{2}\left(  t\right)  ,\,\,\,for\,\,\,t\in\left(
-\infty,0\right]  \cup\left[  \omega_{1},+\infty\right)
\end{array}
\]
It is apparent that $\gamma$ can be approximated by closed geodesics in the
sense of definition \ref{appdef}. We proceed to show that $\gamma$ is not
recurrent by showing that, $\gamma$ and $s\gamma$ are not close in the compact
open topology for any positive real $s.$ For this it suffices to show that
there exists $\varepsilon>0$ and a compact $M\subset\mathbb{R}$ such that for
any positive $s\in\mathbb{R}$,
\begin{equation}
d\bigl(s\gamma\left(  t_{0}\right)  ,\gamma\left(  t_{0}\right)  \bigr)%
\geq\varepsilon\,\,\,for\,\,\,some\,\,\,t_{0}\in M \label{contradiction}%
\end{equation}
For simplicity, we may assume that $d\left(  A,B\right)  =\omega_{1}%
/2=\omega_{2}/4.$ Pick $\varepsilon<d\left(  A,B\right)  /2$ and choose a
compact $M\subset\mathbb{R}$ containing the real numbers $0$ and $3\omega
_{1}/4.$ Let $s$ be arbitrary positive real. If
\[
d\bigl(\gamma\left(  s\right)  ,\gamma\left(  0\right)  \bigr)=d\bigl(%
s\gamma\left(  0\right)  ,\gamma\left(  0\right)  \bigr)\geq\varepsilon
\]
then equation (\ref{contradiction}) is satisfied for the number $t_{0}=0.$ If
$d\bigl(\gamma\left(  s\right)  ,\gamma\left(  0\right)  \bigr)<\varepsilon$
then for $t_{0}=3\omega_{1}/4$ we have
\[
d\bigl(s\gamma\left(  t_{0}\right)  ,\gamma\left(  t_{0}\right)  \bigr)%
=d\bigl(s\gamma\left(  \frac{3\omega_{1}}{4}\right)  ,c_{1}\left(
\frac{3\omega_{1}}{4}\right)  \bigr)>\frac{\omega_{1}}{4}=\frac{d\left(
A,B\right)  }{2}>\varepsilon
\]
This completes the proof that $\gamma$ is not recurrent and, therefore,
approximation by closed geodesics does not imply recurrence.

\medskip

{\footnotesize {\textsc{Agricultural University of Athens, Department of
Mathematics, 75 Iera Odos, Athens 11855 Greece,}\newline \textsc{E-mail:}
\textsl{gmat2xax@auadec.aua.ariadne-t.gr} } }

{\footnotesize {\textsc{University of the Aegean, Department of Mathematics,
Karlovassi, Samos 83200, Greece, E-mail:} \textsl{gtsap@aegean.gr} } }
\end{document}